\documentclass[a4paper,11pt,reqno,twoside]{amsart}

\usepackage{amsmath}
\usepackage{amsfonts}
\usepackage{amssymb}
\usepackage{amsthm}
\usepackage{color}
\usepackage{ifpdf}
\usepackage{array}
\usepackage{url}
\usepackage{multirow,bigdelim}
\usepackage{ascmac}
\usepackage{mathrsfs}
\usepackage[all]{xy} 
\usepackage{graphicx}
\usepackage{ulem}

\newtheorem{theorem}{Theorem}[section]
\newtheorem{proposition}[theorem]{Proposition}
\newtheorem{lemma}[theorem]{Lemma}
\newtheorem{corollary}[theorem]{Corollary}

\newtheorem{remark}[theorem]{Remark}

\numberwithin{equation}{section}

\newcommand*{\C}{\mathbb{C}}
\newcommand*{\R}{\mathbb{R}}
\newcommand*{\Z}{\mathbb{Z}}

\newcommand*{\D}{\mathscr{D}}
\newcommand{\Dom}{\mathfrak D}
\newcommand{\Ker}{\operatorname{Ker}}
\newcommand{\ran}{\operatorname{Ran}}

\newcommand{\PW}{\mathrm{PW}}
\newcommand{\Hglob}{\mathcal H_W}
\newcommand{\Ha}{\mathcal H(A_a)}
\newcommand{\Hhat}{\widehat{\mathcal H(A_a)}}
\newcommand{\Samp}{\mathcal S}
\newcommand{\Jref}{\mathcal R}
\newcommand{\Dlim}{\mathsf D_{\pi/2}}
\newcommand{\Cinv}{\mathsf C}
\newcommand{\Psf}{\Delta_\Gamma}

\newcommand{\comment}[1]{}
\title[]%
      {Spectral convergence associated with \\ Weil's quadratic form} 
\author[M. Suzuki]{Masatoshi Suzuki}
\date{Version of \today}
\AtBeginDocument{%
\begin{abstract}
Assuming the Riemann hypothesis, we study a distinguished
self-adjoint extension of the differentiation operator on
a Hilbert space defined by the Weil quadratic form localized
to a finite interval. 
The local quadratic form is determined by data from finitely many places, 
without requiring the zeros of the Riemann zeta function as input.
We prove that the positive eigenvalues of this extension
converge to the distinct positive
ordinates of the nontrivial zeros of the Riemann zeta function
as the interval expands.
We also prove that the associated normalized spectral entire
functions converge locally uniformly on the complex plane
to the canonical product over these distinct ordinates.
The key step is an identity expressing the inverse of the
local extension as a compression of the inverse of a global
self-adjoint operator.
\end{abstract}
\maketitle
}

\begin{document}

\section{Introduction} \label{section_1}

Let $\zeta(s)$ denote the Riemann zeta function, and set
\[
  \xi(s):=s(s-1)\pi^{-s/2}\Gamma(s/2)\zeta(s),
  \qquad A(z):=\xi(1/2-iz),
\]
where $\xi$ is the Riemann xi-function. Let $\Gamma$ be the set of distinct
zeros of $A$. The Riemann hypothesis (RH) is equivalent to
$\Gamma\subset\R$. Under RH, write
\[
  \Gamma_+:=\Gamma\cap(0,\infty)
  =\{\gamma_1<\gamma_2<\cdots\},
\]
where the zeros are listed without multiplicity. In this paper, we
approximate $\Gamma$ by the spectra of self-adjoint differentiation
operators associated with the Weil quadratic form on finite intervals.
The resulting spectral approximations are determined by local
arithmetic data, without requiring the zeros of $\xi$ as input.

In \cite{Su26}, we considered the unbounded self-adjoint operator $A_a$
on $L^2(-a,a)$ associated with the localized Weil quadratic form introduced
by Connes, Consani, and Moscovici~\cite{CCM25}. Its spectrum is bounded below.
For a real parameter $\lambda$ below the bottom of this spectrum, we
introduced the Hilbert space $\mathcal H(T_a)$ by completing
$C_c^\infty(-a,a)$ with respect to the quadratic form associated with the
positive operator $T_a:=A_a-\lambda I$. 
The differentiation operator $\D_a=i\,d/dx$ on $\mathcal H(T_a)$, with
domain $C_c^\infty(-a,a)$, is then a densely defined symmetric operator
with deficiency indices $(1,1)$.
We studied its family of self-adjoint extensions
$\overline\D_{a,\theta}$, parametrized by a phase $\theta\in[0,2\pi)$.
All these finite-interval results hold without assuming RH. 

Under RH, we may take $\lambda=0$, obtaining the Hilbert space $\Ha$ and
the differentiation operator $\D_a$ on it. 
The first result is as follows.

\begin{theorem} \label{thm_1}
Assume RH, and let $0<t_{a,1}<t_{a,2}<\cdots$ be the positive eigenvalues
of the self-adjoint extension $\overline\D_{a,\pi}$ of $\D_a$ on $\Ha$
corresponding to the phase $\theta=\pi$. Then, for every $n\ge1$,
\begin{equation}
  \label{eq_101}
  t_{a,n}\ge\gamma_n,
  \qquad t_{a,n}\searrow\gamma_n\quad(a\to\infty),
\end{equation}
where $\searrow$ denotes monotone nonincreasing convergence.
\end{theorem}

Each eigenvalue of $\overline\D_{a,\pi}$ is simple, but we do not assume
that the zeros of $\xi$ are simple.

The arithmetic significance of Theorem~\ref{thm_1} is that the local
spectral problems can be defined without prescribing the zeros of
$\xi$. Indeed, for test functions supported in $[-a,a]$, the Weil
explicit formula expresses the localized quadratic form in terms
of the archimedean place and a finite sum over prime powers
$p^m$ satisfying $\log p^m\le 2a$.
Under RH, this quadratic form determines the Hilbert space $\Ha$
and the corresponding self-adjoint differentiation operators.
Although the zeros of $\xi$ enter the spectral representation of
the Weil quadratic form and the proof of convergence, they are not
required as input for the arithmetic definition of the local
spectral problem.

Thus Theorem~\ref{thm_1} provides, for every fixed $n$, a
monotone family of upper spectral approximations to the
distinct zero ordinate $\gamma_n$, obtained by expanding
the support of the localized Weil quadratic form.
The result goes beyond a spectral realization of the set of zeros. 
It establishes a precise convergence relation between
finite-interval spectral problems, determined by local arithmetic
information, and the nontrivial zeros of the Riemann zeta function.
\medskip

The second result is a limit formula for entire functions deduced from
the eigenvalue convergence in Theorem~\ref{thm_1}.
In \cite{Su26}, we introduced an entire function $W(a,\theta,\lambda;z)$
whose zeros are precisely the eigenvalues of $\overline\D_{a,\theta}$.
In this paper, we consider only the case $\lambda=0$
and write $W(a,\theta;z):=W(a,\theta,0;z)$.
For $\theta=\pi$, the function $W(a,\pi;\cdot)$ is even, does not vanish
at the origin, and has as its zeros the eigenvalues of
$\overline\D_{a,\pi}$ (Lemmas~\ref{lem_3_1}
and~\ref{lem_3_2}). Define
\begin{equation} \label{eq_102}
  \Psf(z):=\prod_{\gamma\in\Gamma_+}
  \left(1-\frac{z^2}{\gamma^2}\right)
\end{equation}
by including each distinct positive zero of $A(z)=\xi(1/2-iz)$ once.
The product converges locally uniformly in the whole plane.
Since it does not take into account the multiplicities of the zeros of $\xi$,
the quotient $A/\Psf$ is an entire function that need not be constant.
We have the following theorem.

\begin{theorem} \label{thm_2}
Assume RH. Then
\begin{equation} \label{eq_103}
  \lim_{a\to\infty}\frac{W(a,\pi;z)}{W(a,\pi;0)}=\Psf(z)
\end{equation}
holds uniformly on every compact subset of $\C$.
\end{theorem}

Theorem~\ref{thm_2} extends the convergence of individual eigenvalues 
to convergence at the level of entire functions. 
This is not merely a restatement of Theorem~\ref{thm_1}. 
The bounds $t_{a,n}\ge\gamma_n$ also provide uniform control
of the tails of the corresponding canonical products, allowing
the passage to a locally uniform limit on the whole complex plane.
The spectral function $W(a,\pi;z)$ is constructed from the
reproducing kernel of the Fourier transform of $\Ha$
and is therefore determined by the same local arithmetic
data as the underlying quadratic form.
Consequently, the finite-interval construction yields a family
of entire functions whose normalized limit encodes all the distinct
zero ordinates of $\xi$.
The multiplicities of the zeros are not recovered by this product,
and neither theorem requires the simplicity of the zeta zeros.

Under RH, this result also gives an explicit realization of the
limit formula appearing in \cite[Corollary~1.6, (1.12)]{Su26}.
Namely, there is a family of functions $\phi(a,z)$ holomorphic in the
upper half-plane such that
\begin{equation} \label{eq_104}
  \lim_{a \to \infty}e^{\phi(a,z)}W(a,\pi;z)
  =
  \frac{\xi(1/2-iz)}{\xi(1/2-iz)+\xi'(1/2-iz)}
\end{equation}
holds locally uniformly in the upper half-plane.
The functions $\phi(a,z)$
are given explicitly in Section~\ref{section_5_2}.
Corollary~1.6 of the earlier paper asserts that a limit formula
of this kind, under the admissibility conditions specified there,
implies RH.
Conversely, the present result establishes such a formula under RH
for the distinguished phase $\theta=\pi$ and an explicitly specified
normalization whose dependence on $a$ is only through a nonzero scalar.
Thus we obtain a concrete realization of the corresponding
equivalence with RH.

If RH holds and all nontrivial zeros of $\zeta$ are simple, then
\begin{equation} \label{eq_105}
  \lim_{a\to\infty}\frac{W(a,\pi;z)}{W(a,\pi;0)}
  =\frac{\xi(1/2-iz)}{\xi(1/2)}
\end{equation}
holds locally uniformly on $\C$
(Corollary~\ref{cor_5_2}).
\medskip

The approximation of the nontrivial zeros of $\zeta$
using finite arithmetic information has also been studied
from a different viewpoint by Gonek~\cite{Go12}. 
Using truncated Euler products and the functional equation,
he constructed a family of non-analytic approximations to the Riemann
zeta function and proved that their zeros approach the
zeros of $\zeta$ as the truncation parameter increases under RH.
The present construction differs in that the approximating ordinates
are eigenvalues of self-adjoint differentiation operators.
In addition to their convergence, Theorem~\ref{thm_1} establishes
monotone upper bounds for each distinct zero ordinate, while
Theorem~\ref{thm_2} gives locally uniform convergence of the
associated normalized entire functions to the canonical product
$\Psf$.

Connes and Consani~\cite{CC23} developed a spectral approach based
on the semi-local Weil quadratic form and prolate spheroidal wave
functions. They investigated the small eigenvalues of the localized
form and constructed perturbed Dirac operators whose low-lying
spectra exhibit a close numerical resemblance to the low-lying
zeros of $\zeta$.
Their theory of $\zeta$-cycles also gives exact spectral
correspondences with ordinates of zeros on the critical line
at certain values of the circle length.
Under RH, the present work concerns a different family of self-adjoint
differentiation operators and establishes convergence
to each distinct zero ordinate as the interval expands, 
with a monotone upper bound and 
a corresponding entire-function limit.

A related approach based on the localized Weil quadratic form is
developed by Connes, Consani, and Moscovici~\cite{CCM25}, whose work also
motivated \cite{Su26}. Of particular interest is their proposed connection
between the zeros $\Gamma$ of $\xi$ and the spectrum of a self-adjoint
perturbation of a differentiation operator constructed using $A_a$.
In broad terms, suppressing the technical details of the finite-dimensional
approximation, the argument in \cite{CCM25} shows that, for each $a>0$,
if the smallest eigenvalue
of $A_a$ is simple and its eigenfunction $v_a$ is even, then all zeros of
the Fourier transform $\widehat{v_a}(z)$ are real
(cf.~\cite[Theorem~5.10]{CCM25}).
This is obtained by expressing $\widehat{v_a}(z)$
as the zeta-regularized determinant of the perturbed
operator, up to a nonvanishing exponential factor.
Numerical observations indicate that the zeros of $\widehat{v_a}(z)$
approximate $\Gamma$ well even for small $a>0$.
This suggests that the Fourier transform of the eigenfunction
corresponding to the smallest eigenvalue converges to the Riemann
xi-function after a suitable normalization.
More precisely, one expects a limit formula
\[
  \lim_{a\to\infty}c_a\widehat{v_a}(z)=\xi(1/2-iz)
\]
for suitable normalization constants $c_a$.
Although $W(a,\pi;z)$ and $\widehat{v_a}(z)$ arise from different
constructions, this expected limit is analogous to
\eqref{eq_104} and
\eqref{eq_105}.
The simplicity of the smallest eigenvalue of $A_a$ and the evenness of
the corresponding eigenfunction, assumed in \cite{CCM25}, have been
proved for sufficiently small $a>0$ without assuming RH
(\cite[Theorem~1.4]{Su26}).

We outline the proof strategy using notation introduced in
Section~\ref{section_2}.
Under RH, there is a global Hilbert space $\Hglob$
equipped with a self-adjoint operator $\Dlim$
whose spectrum is $\Gamma$.
This operator can also be realized, via a unitary
isomorphism with a de Branges space $\mathcal H(E)$,
as a self-adjoint extension of the multiplication
operator $F(z)\mapsto zF(z)$.
The essential step is to establish an exact relation
between this global spectral realization and the
finite-interval operators arising from the localized
Weil quadratic form.

Rather than compressing the unbounded operator $\Dlim$ itself,
we compress its inverse.
More precisely, for the natural isometric embedding
$j_a:\Ha\to\Hglob$, the identity
\[
  (\overline\D_{a,\pi})^{-1}
  =j_a^*\Dlim^{-1}j_a
\]
of Theorem~\ref{thm_3_3} is central to the proof of
Theorem~\ref{thm_1}.
The nontrivial point is that the compression of the global inverse
is identified with the inverse of a particular self-adjoint
extension of the local differentiation operator.
With the deficiency vectors fixed by the Riesz representation
in Section~\ref{section_3_1}, the extension is precisely the one
corresponding to $\theta=\pi$.

Writing $\Pi_a=j_aj_a^*$ and $\Cinv=\Dlim^{-1}$, we have
that $\Cinv$ is a Hilbert--Schmidt operator and
$\Pi_a\to I$ in the strong operator topology.
It follows that $\Pi_a\Cinv\Pi_a\to\Cinv$ in Hilbert--Schmidt norm. 
The min--max principle applied to the increasing subspaces $j_a\Ha$
then yields \eqref{eq_101}.
Finally, Hadamard factorization of the even entire function
$W(a,\pi;z)$, together with a uniform estimate for the tails of its
canonical product, gives \eqref{eq_103} and hence
\eqref{eq_104}.
\medskip

The paper is organized as follows. Section~2 fixes the notation and
reviews the constructions of the local and global Hilbert spaces.
Section~3 treats the self-adjoint extensions of the local differentiation
operator and their spectral functions, and
establishes a formula expressing a particular inverse operator
as a compression.
Section~4 uses this identity to prove Theorem~\ref{thm_1}.
Section~5 proves Theorem~\ref{thm_2} and deduces a limit formula
of the type \eqref{eq_104} appearing in
\cite[Corollary~1.6]{Su26}.
Finally, in Section~\ref{section_6}, we prove convergence of finite-dimensional
Riesz approximations to $W(a,\pi;z)$ and report numerical illustrations
based on the arithmetic Gram matrix of the localized Weil form.

\section{Preliminaries} \label{section_2}

Throughout the remainder of the paper, we assume RH unless otherwise stated.

\subsection{Notation and basic facts}
\label{section_2_1}

The $L^2$ inner product $\langle\cdot,\cdot\rangle_{L^2}$ is linear in
the first variable. We write $F^\sharp(z):=\overline{F(\bar z)}$ and
$\C_+:=\{z:\Im z>0\}$. We define the Fourier transform by
\begin{equation} \label{eq_201}
  (\mathcal Fu)(z)=\widehat u(z):=\int_\R u(x)e^{izx}\,dx,
  \qquad \|\mathcal Fu\|_2^2=2\pi\|u\|_2^2.
\end{equation}
On the real axis, the integral is initially defined for
$u\in L^1(\R)\cap L^2(\R)$, and the transform is extended to $L^2(\R)$
by Plancherel's theorem. For compactly supported $L^2$-functions, the
same integral defines an entire function. 
We identify functions on $(-a,a)$ with their extensions
by zero to $\R$, and write
$\PW_a:=\mathcal F L^2(-a,a)$.
For $u\in C_c^\infty(\R)$, integration by parts gives 
$\widehat{u'}(z)=-iz\widehat u(z)$
and $\widehat{iu'}(z)=z\widehat u(z)$.

Under RH, $A$ is a real even entire function. If $m_\gamma$ denotes the
multiplicity of a zero $\gamma\in\Gamma\subset\R$, then
$\Gamma=-\Gamma$ and $m_{-\gamma}=m_\gamma$.
Since $A(0)=\xi(1/2)\ne0$, the fact that $A$ has order $1$ and the
discreteness of its zeros give
\begin{equation} \label{eq_202}
  \sum_{\gamma\in\Gamma}\frac{m_\gamma}{\gamma^2}<\infty.
\end{equation}
For the functional equation and basic properties of the zeros, see
\cite[Sections~10.1--10.2]{MV07}.

\subsection{Local Hilbert spaces and reproducing kernels}
\label{section_2_2}

Unless otherwise stated, the notation and basic facts used here follow
\cite{Su26}. For $u,f\in C_c^\infty(\R)$, 
define the Weil quadratic form by
\begin{equation} \label{eq_203}
  Q_W(u,f):=\sum_{\gamma\in\Gamma}m_\gamma\widehat u(\gamma)
  \overline{\widehat f(\gamma)},\qquad Q_W(u):=Q_W(u,u).
\end{equation}
This is the representation given by the spectral side of the Weil
explicit formula, and $Q_W$ is nonnegative under RH.

The restriction of the Weil quadratic form to $(-a,a)$
is represented by the self-adjoint operator $A_a$
on $L^2(-a,a)$ through
\[
  Q_W(u,f)=\langle A_au,f\rangle_{L^2},
  \qquad u,f\in C_c^\infty(-a,a)
\]
(\cite[Section~3]{CCM25}).
Let $\Ha$ be the $Q_W$-completion of $C_c^\infty(-a,a)$, 
with inner product 
\[
\langle\cdot,\cdot\rangle_{A_a}:=Q_W(\cdot,\cdot).
\]
We use the local coercivity from \cite[Sections~6--7]{Su26} in the form
\begin{equation} \label{eq_204}
  Q_W(u)\ge\lambda_a\|u\|_{L^2(-a,a)}^2,
  \qquad\lambda_a>0,\quad u\in C_c^\infty(-a,a),
\end{equation}
where $\lambda_a:=\inf\sigma(A_a)$. 
By \cite[Section~3.2]{Su26}, $C_c^\infty(-a,a)$ is a core
for the quadratic form associated with $A_a$.
Together with \eqref{eq_204}, this identifies $\Ha$
with the form domain $\Dom(A_a^{1/2})$,
equipped with the $Q_W$-norm. Thus,
\[
  \Ha\hookrightarrow L^2(-a,a),\qquad
  \Hhat:=\mathcal F\Ha\subset\PW_a.
\]
Equip $\widehat{\Ha}$ with the inner product
$\langle\widehat u,\widehat f\rangle_{\widehat A_a}:=\langle u,f\rangle_{A_a}$.
This makes $\widehat{\Ha}$ a de Branges space, as we now show.
For the basic theory of de Branges spaces, see \cite{dB68}; a brief
account is also given in \cite[Section~3.1]{Su23}.

\begin{proposition}
\label{prop_2_1}
The space $\Hhat$ is a de Branges space with no common zeros.
\end{proposition}

\begin{proof}
Since $\Hhat\subset\PW_a$, 
it is a Hilbert space of entire functions. 
We verify the three defining axioms of a de Branges space: 
continuity of point evaluations, invariance under conjugation, 
and invariance under the transformation of nonreal zeros, 
with preservation of the norm in the latter two cases.

For every $z\in\C$,
\[
  |\widehat u(z)|\le
  \left(\int_{-a}^a e^{-2x\Im z}\,dx\right)^{1/2}
  \lambda_a^{-1/2}\|u\|_{A_a},
\]
so point evaluation is continuous.
Conjugation $F\mapsto F^\sharp$ corresponds to
$u(x)\mapsto\overline{u(-x)}$. 
This transformation preserves \(C_c^\infty(-a,a)\) 
and the norm induced by \(Q_W\), 
and hence extends to an isometry on the completion.

Now let $F\in\Hhat$ with $F(w)=0$ and $w\in\C\setminus\R$.
Since $\Hhat\subset\PW_a$, the zero-shifting property of the
Paley--Wiener space gives
\[
  F_w(z):=\frac{z-\bar w}{z-w}F(z)\in\PW_a.
\]
By the characterization of the form domain in
\cite[Section~3]{Su26} and the identification
$\Ha=\Dom(A_a^{1/2})$ above, we have the
set-theoretic equality
\[
  \Ha=
  \left\{u\in L^2(-a,a):
  \int_\R(1+\log^+|t|)
  |\widehat u(t)|^2\,dt<\infty\right\}.
\]
Since $|F_w(t)|=|F(t)|$ on the real axis, this characterization implies
$F_w\in\Hhat$. Moreover, the definition of \(Q_W\) yields 
$\|F_w\|_{\widehat A_a}=\|F\|_{\widehat A_a}$.

Thus all three de Branges axioms hold. For each $z\in\C$, one can choose
a test function whose Fourier transform does not vanish at $z$.
Consequently, the space has no common zeros.
\end{proof}

We follow the indexing convention of \cite[Sections~6.2--6.3]{Su26},
writing $u_{z,a}$ for the vector denoted there by $v_z$ for fixed $a$:
\begin{equation} \label{eq_205}
  \langle u,u_{z,a}\rangle_{A_a}=\widehat u(\bar z)
  \qquad(u\in\Ha).
\end{equation}
The vector $u_{z,a}\in\Ha$ is uniquely determined by
the Riesz representation theorem. 
If $e_z(x):=e^{-izx}$, then
\[
  u_{z,a}=A_a^{-1}e_z\in\Dom(A_a).
\] 
Here and throughout the paper, $\Dom(T)$ denotes
the domain of an operator $T$. 
Indeed, the representation of the closed form gives
\[
  \langle u,A_a^{-1}e_z\rangle_{A_a}
  =
  \langle u,e_z\rangle_{L^2}
  =
  \widehat u(\bar z),
\]
and the assertion follows by uniqueness of the Riesz representative.
The reproducing kernel of the de Branges space $\Hhat$ is
\begin{equation} \label{eq:Ja}
  J_a(w,z):=\langle u_{\bar w,a},u_{\bar z,a}\rangle_{A_a}
  =\widehat{u_{\bar w,a}}(z),\qquad
  J_a(w,w)=\|u_{\bar w,a}\|_{A_a}^2.
\end{equation}
In particular, $J_a(-i,\cdot)=\widehat{u_{i,a}}$ and
$J_a(i,\cdot)=\widehat{u_{-i,a}}$.
Choosing a generating function $E_a$ for $\Hhat$, we have
\begin{equation}
  \label{eq:Ja-Ea}
  J_a(w,z)=\frac{\overline{E_a(w)}E_a(z)-\overline{E_a^\sharp(w)}E_a^\sharp(z)}
  {2\pi i(\bar w-z)}.
\end{equation}

\subsection{The global completion and its model-space realization}
\label{section_2_3}

Let $\Hglob$ be the $Q_W$-completion of $C_c^\infty(\R)$, and denote by
$[u]$ the image of a test function $u$ under the embedding
$C_c^\infty(\R)\hookrightarrow\Hglob$.
Define the weighted sequence space
\[
  m:=(m_\gamma)_{\gamma\in\Gamma},\qquad
  \ell^2(\Gamma,m):=
  \left\{(c_\gamma)_{\gamma\in\Gamma}:
  \sum_{\gamma\in\Gamma}m_\gamma|c_\gamma|^2<\infty\right\}.
\]
The isomorphism we shall use is the well-defined sampling map
\begin{equation} \label{eq:sampling}
  \Samp:\Hglob\overset{\sim}{\longrightarrow}\ell^2(\Gamma,m),\qquad
  \Samp[u]=(\widehat u(\gamma))_{\gamma\in\Gamma}
\end{equation}
from \cite[Lemma~5.4]{Su25}.
The sampling values of a general element of
the completion are defined by taking limits in $\ell^2(\Gamma,m)$.
Consider the natural isometric embedding
\[
  j_a:\Ha~\hookrightarrow~\Hglob,\qquad u\mapsto[u].
\]
Its Hilbert space adjoint
\[
  j_a^*:\Hglob~\to~\Ha
\]
is determined by
\[
  \langle j_a u,h\rangle_{\Hglob}
  =
  \langle u,j_a^*h\rangle_{A_a}\qquad
  (u\in\Ha,\ h\in\Hglob).
\]
Set
\[
  \mathcal M_a:=j_a\Ha.
\]
Since $j_a$ is an isometric embedding, $j_a^*j_a=I_{\Ha}$.
Moreover,
\[
  \Pi_a:=j_aj_a^*
\]
is the orthogonal projection from $\Hglob$ onto $\mathcal M_a$.
Since
\begin{equation} \label{eq_209}
  \mathcal M_a\subset\mathcal M_b\ (a<b),\qquad
  \overline{\bigcup_{a>0}\mathcal M_a}=\Hglob,
\end{equation}
we have $\Pi_a\to I$ in the strong operator topology.
The density follows from
$C_c^\infty(\R)=\bigcup_{a>0}C_c^\infty(-a,a)$.
\medskip

Following \cite{Su26}, set
\[
  E(z):=A(z)+iA'(z)=\xi(1/2-iz)+\xi'(1/2-iz).
\]
Here $A'$ denotes differentiation with respect to $z$, so
$iA'(z)=\xi'(1/2-iz)$. Since RH is assumed, $E$ satisfies the
Hermite--Biehler inequality $|E^\sharp|<|E|$ in the upper half-plane
\cite[Proposition~3.1]{Su23}. 
Let $\mathcal H(E)$ be the de Branges space generated by $E$.
Here we allow the generating function to have real zeros
(cf.~\cite{Su23}). 
The local generating function $E_a$ above must be distinguished from 
the global generating function $E=A+iA'$. 
\medskip

Set $\Theta:=E^\sharp/E$ and consider the model space
\[
  \mathcal K:=K(\Theta):=H^2(\C_+)\ominus\Theta H^2(\C_+)=\mathcal H(E)/E,
\]
where $H^2(\C_+)$ is the Hardy space of the upper half-plane, equipped with
the usual boundary $L^2(\R)$ norm. For basic facts about model spaces,
see \cite[Section~3.1]{Su23}.
An orthonormal basis of $\mathcal K$ is given by
\begin{equation} \label{eq_210}
  F_\gamma(z):=i\sqrt{\frac{m_\gamma}{\pi}}
  \frac{A(z)}{E(z)(z-\gamma)},\qquad
  F_\gamma(\eta)=\frac{\delta_{\gamma\eta}}{\sqrt{\pi m_\gamma}}
\end{equation}
(\cite[Proposition~3.2]{Su23}). Here $\gamma,\eta\in\Gamma$, and the
values at real zeros are defined by removing the singularities.
Following \cite[Section~7.2]{Su26}, define
\begin{equation} \label{eq_211}
\aligned
  U:&\Hglob\overset{\sim}{\longrightarrow}\mathcal K,\qquad
  Uh:=\sum_\gamma\sqrt{m_\gamma}(\Samp h)_\gamma F_\gamma,\qquad \\
  T_E:&\mathcal K\overset{\sim}{\longrightarrow}\mathcal H(E),\quad T_EF=EF.
\endaligned
\end{equation}
Then $U_a:=Uj_a$ is also an isometric embedding.
Let
\[
  V:=\mathcal F^{-1}\mathcal K\subset L^2(0,\infty),
\]
and define $f_\gamma \in V$ by
\[
  \widehat{f_\gamma}=F_\gamma.
\]
An element $h\in\Hglob$ is represented by the function
\begin{equation} \label{eq_212}
  v_h:=\sqrt\pi\,\mathcal F^{-1}Uh
  =\sqrt\pi\sum_\gamma\sqrt{m_\gamma}(\Samp h)_\gamma f_\gamma
\end{equation}
in $V$. We have $\widehat{v_h}(\gamma)=(\Samp h)_\gamma$ and
\[
  \|E\widehat{v_h}\|_{\mathcal H(E)}^2
  =\|\widehat{v_h}\|_2^2=2\pi\|v_h\|_2^2=\pi\|h\|_{\Hglob}^2
\]
by \cite[Theorem~5.5]{Su25}. This representative need not have compact
support, and $U_a u$ is not to be identified with $\widehat u$ as a
function.
\medskip

Define the self-adjoint diagonal operator on $\ell^2(\Gamma,m)$ by
\[
  (T_\Gamma c)_\gamma=\gamma c_\gamma,\qquad
  \Dom(T_\Gamma)=\left\{c:\sum_\gamma m_\gamma\gamma^2|c_\gamma|^2<\infty\right\},
\]
and set
\[
  \Dlim:=\Samp^{-1}T_\Gamma\Samp.
\]
Each distinct point of $\Gamma$ is a simple eigenvalue of this operator.
It corresponds to the self-adjoint extension $M_{\pi/2}$ of the
multiplication operator in \cite[Section~6]{Su25} through
\begin{equation} \label{eq_213}
  (T_EU)\Dlim(T_EU)^{-1}=M_{\pi/2}.
\end{equation}
Indeed, $M_{\pi/2}(EF_\gamma)=\gamma EF_\gamma$, and
$(EF_\gamma)_{\gamma\in\Gamma}$ is an orthonormal basis of
$\mathcal H(E)$; see \cite[Proposition~3.2]{Su23} and
\cite[Section~6]{Su25}.
The eigenfunction expansion in $\mathcal H(E)$ therefore proves
\eqref{eq_213} as an equality of operators, including their domains.
In particular, since $M_{\pi/2}$ is self-adjoint,
$\mathfrak D(M_{\pi/2})$ is dense in $\mathcal H(E)$, and
$\mathfrak D(\Dlim)$ is dense in $\Hglob$. 

Equivalently, one may construct $\Dlim$ by transporting the self-adjoint extension $M_{\pi/2}$ of the multiplication operator on $\mathcal H(E)$ to $\Hglob$ through the unitary isomorphism $T_EU$. Thus, the definition of $\Dlim$ does not intrinsically require prescribing the zero set $\Gamma$ as the spectrum of a diagonal operator. The diagonal realization above is adopted for convenience in the spectral convergence argument.
\medskip

By \eqref{eq_202}, the operator
\[
  \Cinv:=\Dlim^{-1}
\]
is bounded and self-adjoint. Let $\mathfrak S_2(\Hglob)$
denote the Hilbert--Schmidt class on $\Hglob$.
Using the orthonormal basis
$(m_\gamma^{-1/2}\mathbf1_{\{\gamma\}})_{\gamma\in\Gamma}$ of
$\ell^2(\Gamma,m)$, we obtain
\[
  \|\Cinv\|_{\mathfrak S_2}^2
  =\sum_{\gamma\in\Gamma}\frac1{\gamma^2}<\infty.
\]
Thus $\Cinv$ belongs to $\mathfrak S_2(\Hglob)$.
The sum is taken without multiplicity.
Since $0\notin\sigma(\Dlim)$ and the spectral gap
$\inf_{\gamma\in\Gamma}|\gamma|$ is strictly positive,
self-adjointness of $\Dlim$ gives
\[
  \mathfrak D(\Cinv)=\Hglob.
\]

\section{Preparations for the proofs} \label{section_3}

\subsection{The local differentiation operator} 
\label{section_3_1}

Define the differentiation operator $\D_a$ on $\Ha$ by
\[
  \Dom(\D_a):=C_c^\infty(-a,a)\subset\Ha,\qquad \D_a u:=iu'.
\]
The domain $\mathfrak D(\D_a)$ is dense by definition, and $\D_a$ is
symmetric by \cite[Section~6.1]{Su26}. Equation~\eqref{eq_205} gives
\[
  \langle\D_a u,u_{z,a}\rangle_{A_a}
  =\bar z\widehat u(\bar z)=\langle u,zu_{z,a}\rangle_{A_a},
\]
so $\D_a^*u_{z,a}=zu_{z,a}$.
For $f\in\Ker(\D_a^*-z)$, the functional
$L_f(u):=\langle u,f\rangle_{A_a}$ is a distribution on $(-a,a)$.
Indeed, the rapid decay of the Fourier transform of test functions,
together with \eqref{eq_202}, implies that $L_f$ is continuous
in the test-function topology.
Since $iL_f(u')=\bar zL_f(u)$, we have $L_f'=i\bar zL_f$ and hence
$L_f(u)=c\, \widehat{u}(\bar z)$.
Uniqueness of the Riesz representative and 
density of \(C_c^\infty(-a,a)\) in \(\mathcal H(A_a)\) 
now give
\begin{equation} \label{eq_301}
  \Ker(\D_a^*-z)=\C u_{z,a}\quad(z\in\C).
\end{equation}
In particular, $\Ker(\D_a^*\mp i)=\C u_{\pm i,a}$, so the deficiency
indices are $(1,1)$.
\medskip

Reflection 
\[
\Jref u(x):=u(-x)
\] 
is a unitary involution on $\Ha$,
with $\Jref u_{z,a}=u_{-z,a}$ and $\Jref\D_a\Jref=-\D_a$.
Define reflection on the global space by
$(\Samp\Jref f)_\gamma=(\Samp f)_{-\gamma}$. 
Then the injection $j_a$ intertwines the
local and global reflections.
\medskip

The deficiency vectors have equal norms:
\[
  \|u_{i,a}\|_{A_a}^2=J_a(i,i)=J_a(-i,-i)=\|u_{- i,a}\|_{A_a}^2.
\]
Thus von Neumann's theory \cite[Theorem~X.2 and its Corollary]{RSII}
gives all self-adjoint extensions of $\D_a$ as
\begin{equation}
  \label{eq_302}
  \begin{gathered}
  w_{a,\theta}:=u_{i,a}+e^{i\theta}u_{-i,a},\qquad 0\le\theta<2\pi,\\
  \Dom(\overline\D_{a,\theta})
  =\Dom(\overline\D_a)\dotplus\C w_{a,\theta},\qquad
  \overline\D_{a,\theta}:=\D_a^*|_{\Dom(\overline\D_{a,\theta})}, 
  \end{gathered}
\end{equation}
where $\dotplus$ denotes the algebraic direct sum. 
Using the boundary form
\[
  \mathcal W_a(f,h):=\langle\D_a^*f,h\rangle_{A_a}-\langle f,\D_a^*h\rangle_{A_a},
\]
we have 
$f\in\Dom(\overline\D_{a,\theta})$ 
if and only if 
$\mathcal W_a(f,w_{a,\theta})=0$ 
for every $f\in\Dom(\D_a^*)$.

\begin{lemma} \label{lem_3_1}
The domain of $\overline\D_{a,\theta}$ is invariant under reflection
if and only if $\theta=0$ or $\theta=\pi$.
Moreover, $0$ is an eigenvalue of $\overline\D_{a,0}$ but not of
$\overline\D_{a,\pi}$.
\end{lemma}

\begin{proof}
The first assertion follows from
$\Jref w_{a,\theta}=e^{i\theta}w_{a,-\theta}$ and the fact that distinct
phases give distinct domains.
Both $u_{0,a}$ and $w_{a,0}$ are even functions on $(-a,a)$,
whereas $\D_a^*w_{a,0}=i(u_{i,a}-u_{-i,a})$ is odd.
Hence $\mathcal W_a(u_{0,a},w_{a,0})=0$. 

On the other hand, $\Dom(\overline\D_a)\subset H_0^1(-a,a)$ and
$\overline\D_a u=iu'$. Indeed, for a sequence of test functions converging
in the graph norm, applying \eqref{eq_204} to both the functions
and their derivatives gives convergence in $H^1$ and zero boundary values. 
Thus $\overline\D_a$ has no eigenvalues.
Since the domains of two distinct self-adjoint extensions
intersect in $\Dom(\overline\D_a)$, and since
$\Ker(\D_a^*)=\C u_{0,a}$ with
$u_{0,a}\in\Dom(\overline\D_{a,0})$,
the extension $\overline\D_{a,\pi}$ cannot have
$0$ as an eigenvalue.
\end{proof}

\subsection{The spectral function from the boundary form}
\label{section_3_2}

Using the description \eqref{eq_302} of the self-adjoint
extensions, define
\begin{equation}
  \label{eq_303}
  \begin{aligned}
  W(a,\theta;z)
  &:=(z-i)\widehat{u_{i,a}}(z)+e^{i\theta}(z+i)\widehat{u_{-i,a}}(z)\\
  &=(z-i)J_a(-i,z)+e^{i\theta}(z+i)J_a(i,z)\\
  &=\overline{\mathcal W_a(u_{\bar z,a},w_{a,\theta})}.
  \end{aligned}
\end{equation}
The last equality uses our convention that the inner product is linear
in the first variable, together with
\[
  \D_a^*u_{\bar z,a}=\bar z\,u_{\bar z,a},\qquad
  \D_a^*w_{a,\theta}=iu_{i,a}-ie^{i\theta}u_{-i,a},\qquad
  \langle u_{\bar z,a},u_{\pm i,a}\rangle_{A_a}
  =\overline{\widehat{u_{\pm i,a}}(z)}.
\]
This is the spectral function of \cite[(1.11), Section~7.9]{Su26}.

Since $u_{\pm i,a}\in L^2(-a,a)\subset L^1(-a,a)$, their Fourier
transforms are entire, and there is a constant $C_a>0$ such that
\begin{equation}
  \label{eq_304}
  |W(a,\theta;z)|\le C_a(1+|z|)e^{a|\Im z|}
  \qquad(z\in\C,~\theta\in[0,2\pi)).
\end{equation}
Thus $W(a,\theta;\cdot)$ is an entire function of exponential type at
most $a$.

If $W(a,\theta;z)=0$, the boundary form condition implies
$u_{\bar z,a}\in\Dom(\overline\D_{a,\theta})$.
Proposition~\ref{prop_2_1} gives $u_{\bar z,a}\ne0$, so $\bar z$
is an eigenvalue of this self-adjoint operator. Hence $z\in\R$.
Conversely, if $t$ is an eigenvalue of $\overline\D_{a,\theta}$, then
\eqref{eq_301} shows that its eigenfunctions are nonzero multiples
of $u_{t,a}$. Therefore $\mathcal W_a(u_{t,a},w_{a,\theta})=0$, or
equivalently $W(a,\theta;t)=0$.
Thus the zeros of $W(a,\theta;z)$ are precisely the eigenvalues of
$\overline\D_{a,\theta}$.

For the phase $\pi$ used in the proof of Theorem~\ref{thm_1},
we further establish the evenness and simplicity of the zeros of $W(a,\pi;z)$.

\begin{lemma} \label{lem_3_2}
The function $W(a,\pi;z)$ is even and entire, with $W(a,\pi;0)\ne0$.
All its zeros are simple and coincide with the eigenvalues of
$\overline\D_{a,\pi}$.
\end{lemma}

\begin{proof}
The correspondence between the zeros of $W(a,\pi;z)$ and the
eigenvalues of $\overline\D_{a,\pi}$ has already been established.
By reflection symmetry,
\[
  J_a(i,-z)=J_a(-i,z),\qquad
  J_a(-i,-z)=J_a(i,z).
\]
Hence, using \eqref{eq_303}, we obtain
\[
\begin{aligned}
W(a,\pi;-z)
&=(-z-i)J_a(-i,-z)-(-z+i)J_a(i,-z)\\
&=-(z+i)J_a(i,z)+(z-i)J_a(-i,z)=W(a,\pi;z).
\end{aligned}
\]
Thus $W(a,\pi;z)$ is even. 
By Lemma~\ref{lem_3_1}, $0$ is not an eigenvalue of
$\overline\D_{a,\pi}$, and hence $W(a,\pi;0)\ne0$.

It remains to prove that all zeros are simple.
Let $t\in\R$ be a zero of $W(a,\pi;z)$.
Then $u_{t,a}$ is a nonzero eigenvector of
$\overline\D_{a,\pi}$ with eigenvalue $t$.
By the von Neumann decomposition \eqref{eq_302}, write
\[
  u_{t,a}=u_0+c_t w_{a,\pi},
  \qquad u_0\in\Dom(\overline\D_a),\quad c_t\in\C.
\]
Here $c_t\ne0$, since otherwise $u_{t,a}$ would be
an eigenvector of $\overline\D_a$, which has no eigenvalues
by the proof of Lemma~\ref{lem_3_1}.
For every $z\in\C$, the boundary form satisfies
\[
\begin{aligned}
\overline{c_t}\,
\mathcal W_a(u_{\bar z,a},w_{a,\pi})
=\mathcal W_a(u_{\bar z,a},u_{t,a})
=(\bar z-t)
  \langle u_{\bar z,a},u_{t,a}\rangle_{A_a}.
\end{aligned}
\]
Taking complex conjugates and using \eqref{eq_303}
and the defining identity \eqref{eq_205}, we obtain
\[
  W(a,\pi;z)
  =\frac{z-t}{c_t}\widehat{u_{t,a}}(z).
\]
Since $\widehat{u_{t,a}}(t)=\|u_{t,a}\|_{A_a}^2>0$, 
we conclude that
\[
  \left.\frac{\partial}{\partial z}W(a,\pi;z)
  \right|_{z=t}
  =\frac{\|u_{t,a}\|_{A_a}^2}{c_t}\ne0.
\]
Thus every zero of $W(a,\pi;z)$ is simple.
\end{proof}

\subsection{Inverse compression for the distinguished phase}
\label{section_3_3}

We use the isometric embedding $j_a:\Ha\hookrightarrow\Hglob$ 
and the orthogonal projection $\Pi_a=j_aj_a^\ast$ 
introduced in Section~\ref{section_2_3}. 
Recall that $j_a^\ast$ is the Hilbert space adjoint of $j_a$, 
not simply the restriction of a representative function to $(-a,a)$.
\medskip

To prove Theorem~\ref{thm_1}, we compress the inverse
$\Cinv=\Dlim^{-1}$, with $\mathfrak D(\Cinv)=\Hglob$, rather than the
unbounded self-adjoint operator $\Dlim$ itself.
Since $\Cinv$ is bounded and self-adjoint, its compression
$j_a^*\Cinv j_a$ is self-adjoint and defined on all of $\Ha$.
Moreover, as shown in Section~\ref{section_2_3},
$\Cinv$ is Hilbert--Schmidt.
This will yield convergence of the compressed operators
in Hilbert--Schmidt norm. 
We first show that the inverse of the compression is precisely the
self-adjoint extension $\overline{\D}_{a,\pi}$ of the local
differentiation operator $\D_a$.

\begin{theorem} \label{thm_3_3}
With the convention for self-adjoint extensions in
\eqref{eq_302}, we have
\begin{equation} \label{eq_307}
  (\overline\D_{a,\pi})^{-1}=\Cinv_a:=j_a^*\Cinv j_a,\qquad
  \mathfrak D(\Cinv_a)=\Ha.
\end{equation}
The operator $\Cinv_a$ is bounded, self-adjoint, and compact on $\Ha$,
with ${\rm Ker}\,\Cinv_a=\{0\}$ and
$\overline{{\rm Ran}\,\Cinv_a}=\Ha$.
In particular, $0$ belongs to the resolvent set $\rho(\overline\D_{a,\pi})$.
\end{theorem}

\begin{proof}
Since $\Cinv=\Dlim^{-1}$ is self-adjoint and compact and $j_a$ is an
isometry, $\Cinv_a$ is self-adjoint and compact.

We first show that $\Cinv_a\D_a u=u$ for
$u\in C_c^\infty(-a,a)$. For such a function, $iu'$ belongs to the same
test-function space, and
\[
  \sum_\gamma m_\gamma\gamma^2|\widehat u(\gamma)|^2
  =Q_W(iu')<\infty.
\]
By \eqref{eq:sampling} and the definition of $\Dlim$, we therefore have
$j_a u\in\Dom(\Dlim)$ and $\Dlim j_a u=j_a\D_a u$.
It follows that
\begin{equation} \label{eq_308}
  \Cinv_a\D_a u
  =j_a^*\Cinv\Dlim j_a u=j_a^*j_a u=u
  \qquad(u\in C_c^\infty(-a,a)).
\end{equation}

If $f\in\Ker\Cinv_a$, then, for every $u\in C_c^\infty(-a,a)$,
\[
  \langle u,f\rangle_{A_a}
  =\langle\Cinv_a\D_a u,f\rangle_{A_a}
  =\langle\D_a u,\Cinv_a f\rangle_{A_a}=0.
\]
Since $C_c^\infty(-a,a)$ is dense in $\Ha$, we obtain $f=0$. 
Self-adjointness then implies
\[
  (\ran\Cinv_a)^\perp=\Ker\Cinv_a^*=\Ker\Cinv_a=\{0\},
\]
so $\ran\Cinv_a$ is dense in $\Ha$.
Set
\[
  L_a:=\Cinv_a^{-1},\qquad\Dom(L_a):=\ran\Cinv_a.
\]
This is a densely defined self-adjoint operator.
Indeed, symmetry of $L_a$ follows from self-adjointness of $\Cinv_a$.
If $h\in\Dom(L_a^*)$ and $k=L_a^*h$, then, for every $f\in\Ha$,
\[
  \langle f,h\rangle_{A_a}
  =\langle L_a\Cinv_a f,h\rangle_{A_a}
  =\langle\Cinv_a f,k\rangle_{A_a}
  =\langle f,\Cinv_a k\rangle_{A_a}.
\]
Hence $h=\Cinv_a k\in\Dom(L_a)$ and $L_ah=k$, proving $L_a^*=L_a$.
Also, $L_a:\Dom(L_a)\to\Ha$ is bijective and has bounded inverse
$\Cinv_a$, so $0\in\rho(L_a)$.
Equation~\eqref{eq_308} shows both that $u\in\Dom(L_a)$ and
that $L_au=\D_a u$ for every test function $u$.
Thus $L_a$ is a self-adjoint extension of $\D_a$ and hence of its
closure $\overline\D_a$.

It remains to identify the phase $\theta$ for which
$L_a=\overline\D_{a,\theta}$.
For the global reflection, $\Jref\Cinv\Jref=-\Cinv$.
Using $\Jref j_a=j_a\Jref$ and the corresponding adjoint identity gives
$\Jref\Cinv_a\Jref=-\Cinv_a$.
Consequently, $\Dom(L_a)=\ran\Cinv_a$ is invariant under reflection,
and $\Jref L_a\Jref=-L_a$ on this domain.
By Lemma~\ref{lem_3_1}, $L_a$ is the self-adjoint extension
with phase $0$ or $\pi$. The extension with phase $0$ has a zero eigenvalue,
whereas $0\in\rho(L_a)$. Therefore $L_a=\overline\D_{a,\pi}$.
\end{proof}

\section{Proof of Theorem~\ref{thm_1}} \label{section_4}

Consider the operator
\[
  \widetilde\Cinv_a
  :=j_a\Cinv_a j_a^*
  =\Pi_a\Cinv\Pi_a
\]
on the global Hilbert space $\Hglob$.
Since $j_a:\Ha\to\mathcal M_a$ is unitary,
$j_a^*j_a=I_{\Ha}$, and $j_aj_a^*=\Pi_a$, the orthogonal decomposition
$\Hglob=\mathcal M_a\oplus\mathcal M_a^\perp$ gives
\[
  \widetilde\Cinv_a
  \simeq \Cinv_a\oplus 0_{\mathcal M_a^\perp}.
\]
Thus $\widetilde\Cinv_a$ and $\Cinv_a$ have the same nonzero eigenvalues,
with the same multiplicities.
Since $\Pi_a\to I$ in the strong operator topology and
$\Cinv\in\mathfrak S_2$, the standard approximation property of
Hilbert--Schmidt operators gives
\begin{equation} \label{eq_401}
  \|\widetilde\Cinv_a-\Cinv\|_{\mathfrak S_2}
  =\|\Pi_a\Cinv\Pi_a-\Cinv\|_{\mathfrak S_2}
  ~\to~0
  \qquad(a\to\infty).
\end{equation}
In particular, convergence also holds in operator norm.
\medskip

Recall $\Gamma_+=\{\gamma_1<\gamma_2<\cdots\}$.
Since $\Cinv_a$ is injective on the infinite-dimensional space $\Ha$,
it has infinite rank and therefore infinitely many nonzero eigenvalues.
Moreover, $\Cinv_a$ anticommutes with reflection, so its nonzero
eigenvalues occur in positive--negative pairs, with infinitely many of
each sign.
By Theorem~\ref{thm_3_3}, the spectrum of
$\overline\D_{a,\pi}$ consists of the reciprocals of these eigenvalues.
Write its positive eigenvalues as
$0<t_{a,1}<t_{a,2}<\cdots$; its negative eigenvalues are then
$-t_{a,n}$. Each eigenspace is one-dimensional by \eqref{eq_301}.
\medskip

Set $\mu_{a,n}:=1/t_{a,n}$. Then
$\mu_{a,1}>\mu_{a,2}>\cdots>0$ are the positive eigenvalues of
$\Cinv_a$ in decreasing order, and the positive eigenvalues of
$\widetilde\Cinv_a$ form the same sequence.
On the other hand, the positive eigenvalues of $\Cinv=\Dlim^{-1}$ are
$1/\gamma_1>1/\gamma_2>\cdots>0$.

For a compact self-adjoint operator $T$ with at least $n$ positive
eigenvalues, let $\lambda_n^+(T)$ denote the $n$th positive eigenvalue,
listed in decreasing order and counted with multiplicity.
By the max--min principle for compact self-adjoint operators
(see \cite[Theorem~4.22]{Bo20}),
\begin{equation}
  \label{eq_402}
  \lambda_n^+(T)
  =\max_{\substack{L\subset\mathcal H\\\dim L=n}}
  \ \min_{\substack{x\in L\\\|x\|=1}}
  \langle Tx,x\rangle_{\mathcal H},
\end{equation}
where $\mathcal H$ is the Hilbert space on which $T$ acts.
This formula also applies when $T$ has negative eigenvalues.

For $x\in\mathcal M_a$,
\[
  \langle\widetilde\Cinv_a x,x\rangle_{\Hglob}
  =\langle\Pi_a\Cinv\Pi_a x,x\rangle_{\Hglob}
  =\langle\Cinv x,x\rangle_{\Hglob}.
\]
Applying \eqref{eq_402} to the compression on
$\mathcal M_a$ therefore yields
\begin{equation} \label{eq_403}
  \mu_{a,n}
  =\max_{\substack{L\subset\mathcal M_a\\\dim L=n}}
  \ \min_{\substack{x\in L\\\|x\|=1}}
  \langle\Cinv x,x\rangle_{\Hglob}.
\end{equation}
If $0<a<b$, then $\mathcal M_a\subset\mathcal M_b$, so the class of
subspaces over which the maximum is taken increases. Hence
$\mu_{a,n}\le\mu_{b,n}$. Similarly, $\mathcal M_a\subset\Hglob$ gives
\begin{equation} \label{eq_404}
  \mu_{a,n}\le
  \max_{\substack{L\subset\Hglob\\\dim L=n}}
  \ \min_{\substack{x\in L\\\|x\|=1}}
  \langle\Cinv x,x\rangle_{\Hglob}
  =\lambda_n^+(\Cinv)=\frac1{\gamma_n}.
\end{equation}
This proves monotonicity and the upper bound.

To identify the limit, we compare the operators on the common space
$\Hglob$. If compact self-adjoint operators $S$ and $T$ on the same
Hilbert space each have at least $n$ positive eigenvalues, then, for
every unit vector $x$, $|\langle Sx,x\rangle-\langle Tx,x\rangle|\le\|S-T\|$. 
Taking the minimum over unit vectors in each $n$-dimensional subspace,
then the maximum over all such subspaces, gives by \eqref{eq_402}
\begin{equation}
  \label{eq_405}
  |\lambda_n^+(S)-\lambda_n^+(T)|\le\|S-T\|.
\end{equation}
Taking $S=\widetilde\Cinv_a$ and $T=\Cinv$, and using
\eqref{eq_401}, we obtain
\[
  \left|\mu_{a,n}-\frac1{\gamma_n}\right|
  \le\|\widetilde\Cinv_a-\Cinv\|\longrightarrow0.
\]
Thus $\mu_{a,n}\nearrow1/\gamma_n$. Taking reciprocals of these positive
numbers gives $t_{a,n}\ge\gamma_n$ and $t_{a,n}\searrow\gamma_n$.
\hfill $\Box$
\medskip

\begin{remark}
Let $K\subset\C\setminus\Gamma$ be compact.
For each $z\in K$, the operator $I-z\Cinv$ is invertible, and its inverse
is uniformly bounded on $K$.
By \eqref{eq_401} and the Neumann series for the inverse,
$I-z\widetilde\Cinv_a$ is also invertible for all $z\in K$ when $a$ is
sufficiently large, and its inverse converges to $(I-z\Cinv)^{-1}$ in
operator norm, uniformly on $K$.
For such $a$ and $z$, we therefore have
$z\in\rho(\overline\D_{a,\pi})$ and
\[
  j_a(\overline\D_{a,\pi}-z)^{-1}j_a^*
  =\widetilde\Cinv_a(I-z\widetilde\Cinv_a)^{-1}
  \longrightarrow\Cinv(I-z\Cinv)^{-1}=(\Dlim-z)^{-1}
\]
in operator norm, uniformly on $K$.
This compares the local resolvents after transport to $\Hglob$; it is
not the usual resolvent obtained by extending the local operator by
zero on $\mathcal M_a^\perp$. For $z\ne0$, the latter contains the
additional term $-z^{-1}(I-\Pi_a)$.
\end{remark}
\medskip

We now use Theorem~\ref{thm_1} to establish the limit of the associated
spectral entire functions and then derive a limit formula for $A/E$ of
the type in \cite[Corollary~1.6]{Su26}.

\section{Limits of spectral functions} \label{section_5}

\subsection{Proof of Theorem~\ref{thm_2}} 
\label{section_5_1}

By Lemma~\ref{lem_3_2} and
Theorem~\ref{thm_3_3}, $W(a,\pi;\cdot)$ is an even
entire function that does not vanish at the origin, and its zeros are
precisely the simple zeros $\pm t_{a,n}$.
Moreover, Theorem~\ref{thm_1} and \eqref{eq_202} give
\[
  \sum_{n=1}^\infty\frac1{t_{a,n}^2}
  \le\sum_{n=1}^\infty\frac1{\gamma_n^2}<\infty.
\]
Hence the product $\prod_n(1-z^2/t_{a,n}^2)$ over the positive--negative
pairs of zeros converges locally uniformly in the whole plane.

By \eqref{eq_304}, $W(a,\pi;\cdot)$ has order at most $1$.
Applying Hadamard's factorization theorem and grouping the genus-one
canonical factors in positive--negative pairs, we have
\[
  \left(1-\frac z{t_{a,n}}\right)e^{z/t_{a,n}}
  \left(1+\frac z{t_{a,n}}\right)e^{-z/t_{a,n}}
  =1-\frac{z^2}{t_{a,n}^2}.
\]
Thus there are $\alpha_a,\beta_a\in\C$ such that
\[
  W(a,\pi;z)
  =e^{\alpha_a z+\beta_a}
  \prod_{n=1}^\infty\left(1-\frac{z^2}{t_{a,n}^2}\right).
\]
Evenness gives $\alpha_a=0$, and evaluation at $z=0$ yields
$e^{\beta_a}=W(a,\pi;0)$. Consequently,
\begin{equation} \label{eq_501}
  \frac{W(a,\pi;z)}{W(a,\pi;0)}
  =\prod_{n=1}^\infty\left(1-\frac{z^2}{t_{a,n}^2}\right).
\end{equation}

We next justify passing to the limit $a\to\infty$ in this infinite
product. Fix $r>0$, and set
\[
  P_{a,N}(z):=\prod_{n=1}^N\left(1-\frac{z^2}{t_{a,n}^2}\right),\qquad
  P_N(z):=\prod_{n=1}^N\left(1-\frac{z^2}{\gamma_n^2}\right).
\]
By the convergence $t_{a,n}\to\gamma_n$ in Theorem~\ref{thm_1},
$P_{a,N}\to P_N$ uniformly on $|z|\le r$ for each fixed $N$.
The bound $t_{a,n}\ge\gamma_n$ in the same theorem gives
\[
  \sum_{n>N}\frac{|z|^2}{t_{a,n}^2}
  \le\eta_N(r):=r^2\sum_{n>N}\frac1{\gamma_n^2}\longrightarrow0
  \qquad(N\to\infty),
\]
uniformly for $a>0$ and $|z|\le r$.
To estimate the tails of the products themselves, we use the inequality
\[
  \left|\prod_k(1+b_k)-1\right|
  \le\prod_k(1+|b_k|)-1
  \le\exp\!\left(\sum_k|b_k|\right)-1,
\]
valid whenever $\sum_k|b_k|<\infty$.
It follows that, for $|z|\le r$,
\[
  \left|\prod_{n>N}\left(1-\frac{z^2}{t_{a,n}^2}\right)-1\right|
  \le e^{\eta_N(r)}-1.
\]
The same estimate holds with $t_{a,n}$ replaced by $\gamma_n$.
We also have the bound
\[
  |P_{a,N}(z)|,\ |P_N(z)|
  \le B_r:=\exp\!\left(r^2\sum_{n=1}^\infty\frac1{\gamma_n^2}\right),
\]
uniformly in $a$ and $N$.
Separating each product into its first $N$ factors and its tail gives
\[
\begin{aligned}
  \sup_{|z|\le r}
  \left|\frac{W(a,\pi;z)}{W(a,\pi;0)}-\Psf(z)\right|
  &\le 2B_r\bigl(e^{\eta_N(r)}-1\bigr)\\
  &\quad+\sup_{|z|\le r}|P_{a,N}(z)-P_N(z)|.
\end{aligned}
\]
First choosing $N$ sufficiently large and then letting $a\to\infty$
makes the right-hand side arbitrarily small. Since $r$ is arbitrary,
convergence is locally uniform in the whole plane, including on compact
sets that contain zeros.
\hfill $\Box$
\medskip

In this proof, convergence of each eigenvalue gives convergence of the
finite products, while $t_{a,n}\ge\gamma_n$ provides an estimate for
the tails that is uniform in $a$.
Multiplying the two deficiency vectors $u_{i,a}$ and $u_{-i,a}$ in
${\rm Ker}\,(\D_a^\ast\mp i)$ by the same nonzero constant does not
change the quotient $W(a,\pi;z)/W(a,\pi;0)$.

\subsection{A limit formula of the form \eqref{eq_104}}
\label{section_5_2}

As a consequence of Theorem~\ref{thm_2}, we obtain a limit formula
of the form \eqref{eq_104}, which appears in
\cite[Corollary~1.6]{Su26}.
Grouping the Hadamard product of $A$ into positive--negative pairs of
zeros, and using evenness of $A$ and normalization at the origin, gives
the factorization
\begin{equation} \label{eq_502}
  \begin{gathered}
  G(z):=\prod_{\gamma\in\Gamma_+}
  \left(1-\frac{z^2}{\gamma^2}\right)^{m_\gamma-1},\qquad
  A_0:=\frac{A}{G}=A(0)\Psf,\qquad E_0:=\frac{E}{G}.
  \end{gathered}
\end{equation}
In particular, $A=GA_0$ and $E=GE_0$.
Since the order of each zero of $G$ is $m_\gamma-1$, the function $E_0$
is entire and has no zeros on the real axis.
As $E$ satisfies the Hermite--Biehler inequality in $\C_+$, $E_0$ is a
Hermite--Biehler function without real zeros, and
$E_0^\sharp/E_0=E^\sharp/E=\Theta$.
The space $\mathcal H(E)$ used above, allowing real zeros of the
generating function, is understood as $G\mathcal H(E_0)$.

\begin{corollary} \label{cor_5_1}
Assume RH, and let $E_0$ be the function defined in
\eqref{eq_502}. Then
\begin{equation} \label{eq_503}
  \lim_{a\to\infty}\frac{\xi(1/2)}{W(a,\pi;0)E_0(z)}W(a,\pi;z)
  =
  \frac{\xi(1/2-iz)}{\xi(1/2-iz)+\xi'(1/2-iz)}
\end{equation}
uniformly on every compact subset of the upper half-plane.
\end{corollary}

\begin{proof}
We know that $W(a,\pi;0)\ne0$. Using Theorem~\ref{thm_2} and the
identities $A(0)\Psf=A/G$ and $E=GE_0$, we obtain
\[
  \frac{A(0)}{E_0(z)}\frac{W(a,\pi;z)}{W(a,\pi;0)}
  \longrightarrow\frac{A(0)\Psf(z)}{E_0(z)}=\frac{A(z)}{E(z)}.
\]
Since $E_0$ does not vanish in the upper half-plane, convergence is
uniform on every compact subset of that half-plane.
\end{proof}

The normalization factor on the left-hand side of
\eqref{eq_503} is the product of the scalar
$1/W(a,\pi;0)$, which depends only on $a$, and the function $A(0)/E_0(z)$,
which does not.
Fix a holomorphic logarithm of $E_0$ in the upper half-plane, and choose
a logarithm of the scalar for each $a$. Then
\begin{equation} \label{eq_504}
  \phi(a,z):=\log\!\left(\frac{\xi(1/2)}{W(a,\pi;0)}\right)-\log E_0(z)
\end{equation}
is holomorphic in $z$, and \eqref{eq_503} takes the
same form as \eqref{eq_104}, with an exponential
factor. Thus, under RH, we obtain a limit formula of the type in
\cite[Corollary~1.6, (1.12)]{Su26}. 
In the notation of \cite[Corollary~1.6]{Su26},
\eqref{eq_503} corresponds to the choices
$\lambda(a)=0$, $v_\pm(a,\cdot)=u_{\pm i,a}$,
and $\theta(a)=\pi$.

If all zeros of $\xi(s)$ are simple, then $G=1$, and
\eqref{eq_502} and Theorem~\ref{thm_2} give the following.

\begin{corollary} \label{cor_5_2}
Assume RH and that all zeros of $\xi$ are simple. 
Then \eqref{eq_105} holds locally uniformly on $\C$.
\end{corollary}
\medskip

If arbitrary holomorphic functions $\phi(a,z)$ are allowed in
\eqref{eq_104}, then under RH one can obtain an
identity simply by taking a logarithm of the quotient of $A/E$ and
$W(a,\pi;\cdot)$, since both are nonvanishing in the upper half-plane.
Corollaries~\ref{cor_5_1} and~\ref{cor_5_2} do not use this
procedure. Instead, we first prove \eqref{eq_103} and
then restrict the $a$-dependence in \eqref{eq_504} to a single
scalar.

\section{Numerical experiments}
\label{section_6}

Following the finite-dimensional approach of
\cite[Sections~3--6]{CCM25}, we restrict the localized Weil form to
trigonometric polynomials. Instead of the Fourier transform of a
lowest eigenvector of its Gram matrix, we use Riesz vectors to
approximate the spectral function \eqref{eq_303}.
Assume RH and fix $a>0$. For $N\in\Z_{\ge0}$, set
\begin{equation}\label{eq:num-basis}
  \chi_{n,a}(x):=\frac{e^{i\pi nx/a}}{\sqrt{2a}}
  \mathbf1_{[-a,a]}(x),\qquad
  \mathcal X_{a,N}:=\operatorname{span}\{\chi_{n,a}:n\in\Z,\ |n|\le N\}.
\end{equation}
These functions belong to $\Ha$ by its description in the proof of Proposition~\ref{prop_2_1}; 
no membership in $\Dom(\D_a)$ is required.
Write
\begin{equation}\label{eq:num-gram}
\begin{aligned}
  G_{a,N}
  &:=\bigl(Q_W(\chi_{m,a},\chi_{n,a})\bigr)_{m,n=-N}^{N},\\
  \varphi_{n,a}(z)
  &:=\widehat{\chi_{n,a}}(z)
  =\sqrt{2a}\,\frac{\sin(az+n\pi)}{az+n\pi}.
\end{aligned}
\end{equation}
The positive definite matrix $G_{a,N}$ is calculated from
\cite[(5.15) and (5.16)]{Yo92}, specialized to $\mathbb Q$
($r_1=1$, $r_2=0$, and $\log A_{\mathbb Q}=-\log\pi$
in the notation of that paper). 
Besides the pole and archimedean terms, only prime powers
$p^r\le e^{2a}$ occur: the correlations
$\chi_{m,a}*\widetilde{\chi_{n,a}}$ are supported in $[-2a,2a]$,
where $\widetilde f(x)=\overline{f(-x)}$.
They vanish at $\pm2a$, so equality at the cutoff contributes zero.
The computation of $G_{a,N}$ does not require RH
or any information about the zeros of $\xi(s)$.

Let $u^{(N)}_{z,a}\in\mathcal X_{a,N}$ satisfy
$Q_W(v,u^{(N)}_{z,a})=\widehat v(\bar z)$ for
$v\in\mathcal X_{a,N}$.
Writing 
\[
u^{(N)}_{z,a}=\sum_{n=-N}^{N}c_n(z)\chi_{n,a},\]
set
\[
  \boldsymbol c(z)=(c_n(z))_{n=-N}^{N},\qquad
  \boldsymbol\varphi_{a,N}(z)=(\varphi_{n,a}(z))_{n=-N}^{N},
\]
as column vectors. Then
\begin{equation}\label{eq:num-riesz}
  G_{a,N}\overline{\boldsymbol c(z)}
  =\boldsymbol\varphi_{a,N}(\bar z).
\end{equation}
The conjugation follows from our convention that the inner product is linear in the first variable. 
Define
\begin{equation}\label{eq:num-kernel-W}
\begin{aligned}
  J_{a,N}(w,z)&:=\widehat{u^{(N)}_{\bar w,a}}(z)
  =\boldsymbol\varphi_{a,N}(z)^{\mathsf T}
     \overline{G_{a,N}^{-1}\boldsymbol\varphi_{a,N}(w)},\\
  W_N(a,\pi;z)&:=(z-i)J_{a,N}(-i,z)-(z+i)J_{a,N}(i,z).
\end{aligned}
\end{equation}
Thus, $W_N$ can be evaluated by solving \eqref{eq:num-riesz} for $z=\pm i$. 

\begin{proposition}\label{prop:num-convergence}
Assume RH and fix $a>0$. Then $u^{(N)}_{z,a}\to u_{z,a}$ in $\Ha$,
uniformly for $z$ in compact subsets of $\C$.
Moreover, $J_{a,N}\to J_a$ uniformly on compact subsets of $\C^2$,
and $W_N(a,\pi;\cdot)\to W(a,\pi;\cdot)$ locally uniformly on $\C$.
For all sufficiently large $N$, $W_N(a,\pi;0)\ne0$, and
\begin{equation}\label{eq:num-normalized-limit}
  \lim_{N \to \infty}\frac{W_N(a,\pi;z)}{W_N(a,\pi;0)}
  = \frac{W(a,\pi;z)}{W(a,\pi;0)}
\end{equation}
holds uniformly on every compact subset of $\mathbb C$.
\end{proposition}
\begin{proof} 
By \cite[Lemma~2.2]{CC23},
the linear span of $\{e^{i\pi nx/a}:n\in\mathbb Z\}$
is a core for the quadratic form $Q_W^{\,a}$.
By \cite[Section~3.2, proof of Theorem~1]{Su26},
$C_c^\infty(-a,a)$ is also a core for this form.
Hence, under RH, \eqref{eq_204} and the definition of
$\Ha$ imply
\[
  \overline{\bigcup_{N\ge0}\mathcal X_{a,N}}^{\,\Ha}
  =\Ha.
\]
Let $\mathsf P_{a,N}$ be the $Q_W$-orthogonal projection onto
$\mathcal X_{a,N}$. Uniqueness in the Riesz representation gives
\[
  u^{(N)}_{z,a}=\mathsf P_{a,N}u_{z,a},\qquad
  \mathsf P_{a,N}~\to~ I\quad\hbox{strongly}.
\]
For $e_z(x)=e^{-izx}$, \eqref{eq_204} and \eqref{eq_205} imply
\[
  \|u_{z,a}-u_{w,a}\|_{A_a}
  \le\lambda_a^{-1/2}\|e_z-e_w\|_{L^2(-a,a)}.
\]
Hence, the strong convergence of $\mathsf P_{a,N}$
is uniform on $\{u_{z,a}:z\in K\}$ for every
compact subset $K$ of $\C$. Orthogonality also gives
\[
  J_a(w,z)-J_{a,N}(w,z)
  =\left\langle(I-\mathsf P_{a,N})u_{\bar w,a},
                    (I-\mathsf P_{a,N})u_{\bar z,a}\right\rangle_{A_a},
\]
which proves the assertion for $J_{a,N}$ by Cauchy--Schwarz.
The assertions for $W_N$ and its normalization follow from 
\eqref{eq_303}, \eqref{eq:num-kernel-W}, and Lemma~\ref{lem_3_2}.
\end{proof}

We computed $G_{a,N}$ and $W_N$ in Mathematica~11 for
\[
  a=\frac12\log X,\qquad X=3,5,7,11,\qquad N=14,16,18.
\]
Only prime powers $p^r<X$ contribute to the Gram matrices,
since the contributions at $p^r=X$ vanish.
The Gram matrices were evaluated at 140 and 180 decimal-digit
precision, with the archimedean series truncated so that
the entrywise tail bound is less than $10^{-188}$.
Since $-iW_N(a,\pi;t)$ is real for real $t$, we searched
for its zeros by detecting sign changes on a grid of
4001 points in $[0,40]$ and then applying \texttt{FindRoot}
at 140-digit precision.
Table~\ref{tab:num-zeros} compares the first six detected positive
zeros of $W_{18}(a,\pi;\cdot)$ with those of $A$.

\begin{table}[htbp]
\caption{First six positive zero candidates of $W_{18}(\tfrac12\log X,\pi;z)$,
compared with the zeros $\gamma_j$ of $A(z)$ (10 decimal places).}
\label{tab:num-zeros}
\centering
\scriptsize
\setlength{\tabcolsep}{3pt}
\begin{tabular}{crrrrr}
\hline
$j$ & $\gamma_j$ & $X=3$ & $X=5$ & $X=7$ & $X=11$\\
\hline
1 & 14.1347251417 & 14.1347274566 & 14.1347251417 & 14.1347251417 & 14.1347251417\\
2 & 21.0220396388 & 21.0226378897 & 21.0220396388 & 21.0220396388 & 21.0220396388\\
3 & 25.0108575801 & 25.0192189912 & 25.0108575802 & 25.0108575801 & 25.0108575801\\
4 & 30.4248761259 & 30.6272491287 & 30.4248761361 & 30.4248761259 & 30.4248761259\\
5 & 32.9350615877 & 33.6051319793 & 32.9350616846 & 32.9350615877 & 32.9350615877\\
6 & 37.5861781588 & 38.6158405572 & 37.5861810892 & 37.5861781588 & 37.5861781588\\
\hline
\end{tabular}
\end{table}

Let $t_{a,N,j}$ denote the $j$-th positive zero of
$W_N(a,\pi;\cdot)$ found numerically, and put
\[
  d_6^{(N)}=t_{a,N,6}-\gamma_6.
\]
To measure the sensitivity of the linear system
\eqref{eq:num-riesz} to perturbations, we use the condition number 
$
  \kappa_2(G_{a,N})
  :=\lambda_{\max}(G_{a,N})/\lambda_{\min}(G_{a,N})$.
Table~\ref{tab:num-errors} gives the differences
$d_6^{(N)}$ computed from the 20-decimal-place values,
together with the condition numbers.

\begin{table}[htbp]
\caption{Differences $d_6^{(N)}$ for $a=\frac12\log X$
and condition numbers $\kappa_2(G_{a,18})$.}
\label{tab:num-errors}
\centering
\scriptsize
\setlength{\tabcolsep}{5pt}
\begin{tabular}{crrrr}
\hline
$X$ & $d_6^{(14)}$ & $d_6^{(16)}$ & $d_6^{(18)}$
& $\kappa_2(G_{a,18})$\\
\hline
3  & $1.04$ & $1.03$ & $1.03$ & $4.70\times10^{7}$\\
5  & $5.21\times10^{-6}$ & $3.00\times10^{-6}$
   & $2.93\times10^{-6}$ & $1.85\times10^{17}$\\
7  & $1.11\times10^{-9}$ & $1.52\times10^{-11}$
   & $5.43\times10^{-13}$ & $8.24\times10^{25}$\\
11 & $1.28\times10^{-9}$ & $3.70\times10^{-15}$
   & $1.03\times10^{-17}$ & $2.33\times10^{34}$\\
\hline
\end{tabular}
\end{table}

For $N=18$, the differences decrease rapidly as $X$ increases.
At $X=3$, however, the sixth zero still changes by about
$1.93\times10^{-3}$ between $N=16$ and $18$,
although this change is not visible at the precision
shown in Table~\ref{tab:num-errors}.
The condition numbers also increase considerably with $X$.
Thus the numerical agreement should not be interpreted as
a rigorous error bound.

Although not included in the tables, an additional zero
of $W_{14}$ was found at
$z=30\pi/\log11\simeq 39.3044$ for $X=11$.
This zero arises from
$\varphi_{n,a}(z)=0$ for all $|n|\le14$.

As in \cite{CCM25}, these computations show that
low-lying zeros of $A$ can be detected with high accuracy
from finite-dimensional Weil forms.
Proposition~\ref{prop:num-convergence} establishes convergence
to $W(a,\pi;\cdot)$ for fixed $a$, while
Theorems~\ref{thm_1} and~\ref{thm_2} concern the limit
$a\to\infty$.

\bigskip
\noindent\textbf{Acknowledgment.}
This work was supported by JSPS KAKENHI Grant Number JP23K03050.

\bigskip

\noindent
Department of Mathematics, \\
School of Science, \\ 
Institute of Science Tokyo \\
2-12-1 Ookayama, Meguro-ku, \\
Tokyo 152-8551, Japan  \\[2pt]
Email: {\tt msuzuki@math.sci.isct.ac.jp}


\begin{thebibliography}{99}
%
\bibitem{Bo20}
D. Borthwick,
\newblock{Spectral Theory: Basic Concepts and Applications},
Graduate Texts in Mathematics, vol.~284,
Springer, Cham, 2020.
%
\bibitem{CC23}
A. Connes, C. Consani,
\newblock{Spectral triples and $\zeta$-cycles},
Enseign. Math. \textbf{69} (2023),
no.~1--2, 93--148.
%
\bibitem{CCM25}
A. Connes, C. Consani, H. Moscovici,
\newblock{Zeta spectral triples},
in \newblock{Applications of Noncommutative Geometry to Gauge Theories,
Field Theories, and Quantum Space-Time},
EMS Ser. Lect. Math. \textbf{37} (2026).
%
\bibitem{dB68}
L. de Branges,
\newblock{Hilbert Spaces of Entire Functions},
Prentice-Hall, Englewood Cliffs, N.J., 1968.
%
\bibitem{Go12}
S. M. Gonek,
\newblock{Finite Euler products and the Riemann hypothesis},
Trans. Amer. Math. Soc. \textbf{364} (2012),
no.~4, 2157--2191.
%
\bibitem{MV07}
H. L. Montgomery, R. C. Vaughan, 
\newblock{Multiplicative number theory. I. Classical theory, 
Cambridge Studies in Advanced Mathematics, 97
}, 
\newblock{\it Cambridge University Press, Cambridge}, 
\newblock{2007}. 
%
\bibitem{RSII}
M. Reed, B. Simon,
\newblock{Methods of Modern Mathematical Physics II:
Fourier Analysis, Self-Adjointness},
Academic Press, New York, 1975.
%
\bibitem{Su23}
M. Suzuki, \newblock{Li coefficients as norms of functions in a model space},
J. Number Theory \textbf{252} (2023), 177--194.
%
\bibitem{Su25}
M. Suzuki, \newblock{On the Hilbert space derived from the Weil distribution},
Canadian Journal of Mathematics, First View.
%
\bibitem{Su26}
M. Suzuki, \newblock{Weil's quadratic form via the screw function},
\url{https://arxiv.org/abs/2606.09096}.
%
\bibitem{Yo92}
H. Yoshida,
\newblock{On Hermitian forms attached to zeta functions},
in \newblock{Zeta Functions in Geometry},
Adv. Stud. Pure Math. \textbf{21} (1992), 281--325.
%
\end{thebibliography}
\end{document}